\documentclass[12pt]{article}
\usepackage{amssymb,amsfonts,amsmath}
\usepackage{mathrsfs}

\def\pf{\noindent {\it Proof.} }
\def\qed{\hfill $\square$ }

\def\R{{\mathbb R}}
\def\C{\mathcal{C}}
\def\L{\mathfrak{L}}

\makeatletter\@addtoreset{figure}{section}\makeatother

\newcommand{\captionfonts}{\small}
\makeatletter
\long\def\@makecaption#1#2{%
   \vskip 10\p@
   \setbox\@tempboxa\hbox{\captionfonts {#1} \ \ #2}%
   \ifdim \wd\@tempboxa >\hsize
       \captionfonts {#1}\ \ #2\par
   \else
       \hbox to\hsize{\hfil\box\@tempboxa\hfil}%
   \fi}
\makeatother

\newtheorem{theo}{Theorem}

\newtheorem{coro}{Corollary}[section]

\newtheorem{lem}{Lemma}[section]

\makeatletter \@addtoreset{equation}{section}
\@addtoreset{theo}{section} \makeatother

\allowdisplaybreaks

\title{Koroljuk's formula for counting lattice paths revisited}
\author{James J.Y. Zhao\\
\small Dongling School of Economics and Management,\\
\small University of Science and Technology Beijing,\\
\small Beijing 100083, P.R. China\\
\small and\\
\small Center for Applied Mathematics,\\
\small Tianjin University, \\
\small Tianjin 300072, P.R. China\\
\small \texttt{Email:\,zhaojy@ustb.edu.cn}}

\date{\small{Jun 25, 2013}\\
Mathematics Subject Classification: 05A15, 05A19}

\begin{document}

\maketitle

\begin{abstract}
Koroljuk gave a summation formula for counting the number of lattice paths from $(0,0)$ to $(m,n)$ with $(1,0), (0,1)$-steps in the plane that stay strictly above the line $y=k(x-d)$, where $k$ and $d$ are positive integers. In this paper we obtain an explicit formula for the number of lattice paths from $(a,b)$ to $(m,n)$
above the diagonal $y=kx-r$, where $r$ is a rational number. Our result slightly generalizes  Koroljuk's formula, while the former can be essentially derived from the latter. However, our proof uses a recurrence with respect to the starting points, and hereby presents a new approach to Koroljuk's formula.
\end{abstract}

{\bf Keywords:} Koroljuk's formula, forward recursion, lattice path, ballot problem, reflection principle, Fuss-Catalan Number, Hagen-Rothe identity.

\section{Introduction}

Lattice paths are simple and elementary combinatorial structures, while they have many interesting applications in combinatorics, probability and statistics, see for instance Mohanty \cite{Mohanty79} and Krattenthaler and Mohanty \cite{krat-Moha}.
In the study of the discrepancy of empiric distributions, Koroljuk \cite{Koroljuk} obtained a formula for counting the number of lattice paths in the plane that stay strictly above certain line with slope and $y$-intercept corresponding to a pair of integers. The main objective of this paper is to give a new approach to Koroljuk's formula.
The original proof given by Koroljuk used a technique of decomposing the paths, and our proof relies on a recurrence with respect to the starting points.

Let us first give an overview of Koroljuk's formula. It should be mentioned that Koroljuk's formula is originally concerned with the enumeration of lattice paths from $(0,0)$ to $(m-pn,m+n)$ using the steps $(1,1)$ and $(-p,1)$ in the plane that intersect the line $x=c$, for some positive integers $m,n,p$ and $c$. Koroljuk \cite[eq. (4)]{Koroljuk} \footnote{Here the original upper bound $m+n-[\frac{c}{p}]$ of summation is changed to  $c+\lfloor\frac{m+n-c}{p+1}\rfloor(p+1)$ without affecting the value of the summation.} found that such lattice paths are counted by
\begin{equation}\label{eq-Koro}
\sum\limits_{s\equiv c \!\!\! \pmod{\!(p+1)}}^{c+\lfloor\frac{m+n-c}{p+1}\rfloor(p+1)}\frac{c}{s}{s\choose \frac{s-c}{p+1}}{m+n-s\choose n-\frac{s-c}{p+1}},
\end{equation}
where $\lfloor x\rfloor$ denotes the largest integer not greater than $x$.

Koroljuk's formula also counts the number of certain lattice paths using $(1,0), (0,1)$-steps in the plane with some restrictions, see for instance Niederhausen \cite[p. 9]{Nieder02}.  For some positive integers $m, k$ and rational number $d\geq \frac{k-1}{k}m$, Niederhausen \cite[p. 9]{Nieder02} \footnote{Here the original lower bound $\lfloor d\rfloor+1$ of summation is changed to $\lfloor\frac{kd-1}{k+1}\rfloor +1$, which is more accurate. The details are shown in Section \ref{sec-re}.} showed that the number of lattice paths from $(0,0)$ to $(m,n)$ with $(1,0), (0,1)$-steps that stay strictly  above the line $y=k(x-d)$ is equal to
\begin{equation}\label{eq-NCp9}
{m+n\choose m}-\sum\limits_{i=\lfloor\frac{kd-1}{k+1}\rfloor +1}^m \frac{n-k(m-d)}{n-k(i-d)}{i+k(i-d)\choose i}{m-i-1+n-k(i-d)\choose m-i}.
\end{equation}
Substituting $p=k$, $m=n$, $n=m$, $c=n-k(m-d)$ and $s=n-k(m-d)+(k+1)i$
in Koroljuk's formula \eqref{eq-Koro}, we obtain
\begin{align*}
\sum\limits_{i=0}^{m-\lceil \frac{kd}{k+1}\rceil}\frac{n-k(m-d)}{n-k(m-d)+(k+1)i}
{n-k(m-d)+(k+1)i\choose i}{m+k(m-d)-(k+1)i\choose m-i},
\end{align*}
where $\lceil x\rceil$ is the least integer not less than $x$.
By further substituting $m-i$ for $i$, the above formula becomes
\begin{align*}
\sum\limits_{i=\lceil \frac{kd}{k+1}\rceil}^m \frac{n-k(m-d)}{m+n+kd-(k+1)i}
{m+n+kd-(k+1)i\choose m-i}{(k+1)i-kd\choose i}.
\end{align*}
It is easy to verify that the above sum is equal to the subtrahend term of \eqref{eq-NCp9}. Therefore, Koroljuk's formula also gives
the number of lattice paths from $(0,0)$ to $(m,n)$ with $(1,0), (0,1)$-steps that touch or cross the line $y=k(x-d)$. The equivalence between these two kinds of lattice paths counted by Koroljuk's formula will be illustrated in the proof of Theorem \ref{theo-3.1}.

In this paper, we are concerned with the enumeration of lattice paths from $(a,b)$ to $(m,n)$, using $(1,0), (0,1)$-steps in the plane, that stay above or strictly above a boundary line $y=kx-r$, where $a,b,m,n,k$ are nonnegative integers and $r$ is a rational number. Although the enumeration of lattice paths starting from a non-origin point can be transformed to that from the origin, it is convenient to directly work with the former. Based on this consideration, we present another approach to Koroljuk's formula. The main feature of the present approach is to use the forward recursion with respect to the starting points of lattice paths. By using some reflections on lattice paths, we derive the enumeration formula for the lattice paths with respect to the boundary line $y=x/k-r$ from that with respect to $y=kx-r$.

We would like to point out that the enumeration problem considered here has attracted much attention over the past several decades. We shall give an overview of some main methods which have been used in this field, and these methods vary significantly depending upon the boundary lines and the starting points.
For the convenience, let us first introduce some notations. Given rational numbers $q>0$ and $r$, and integers $a,b,m,n$ with $0\leq a\leq m$ and $0\leq b\leq n$, let $L_{q,r}(a,b;m,n)$ denote the set of lattice paths from $(a,b)$ to $(m,n)$ that stay above the line  $y=qx-r$, and let $\L_{q,r}(a,b;m,n)$ denote the set of lattice paths from $(a,b)$ to $(m,n)$ that stay strictly above the line $y=qx-r$. With these notations,  Koroljuk's formula is related to the evaluation of $|\L_{k,kd}(0,0; m,n)|$.

The first lattice path enumeration problem is to count $|L_{1,0}(0,0; m,n)|$. It is often thought that this problem  originated from Bertrand's ballot problem \cite{Bertrand}, with the solution
$$|L_{1,0}(0,0; m,n)|=\frac{n-m+1}{n+1}{m+n\choose m},$$
which, known as the ballot number, is a generalization of the well known Catalan number $\frac{1}{n+1}{2n\choose n}$.
One of the most classical techniques for solving Bertrand's ballot problem is the renowned \emph{reflection principle}, which has been widely contributed to Andr\'{e} \cite{Andre1887}. Renault \cite{Renault} pointed out that the reflection principle is in fact a variation due to Aebly \cite{Aebly1923} and Mirimanoff \cite{Miri1923}, see also \cite[\S 2.3]{Hump10}.

Barbier \cite{Barbier} further studied the generalized ballot problem, which is equivalent to that of counting $|L_{k,0}(0,0; m,n)|$.
The solution
\begin{equation}\label{eq-GBP}
|L_{k,0}(0,0; m,n)|={m+n\choose m}-k{m+n\choose m-1}
\end{equation}
is known as the generalized ballot number. In the spirit of the reflection principle, Goulden and Serrano \cite{Gou-Ser} gave a bijective proof for the generalized ballot problem by using rotation instead of reflection.

It is clear that the boundary line for the ballot problem or its generalization has zero $y$-intercept. If the $y$-intercept of the boundary line is nonzero, the corresponding lattice path enumeration problem will become a little complicated. As remarked before, Koroljuk \cite{Koroljuk} essentially studied the problem of counting $|\L_{k,kd}(0,0; m,n)|$, where $k,d$ are positive integers. Niederhausen \cite[Example 2.2]{Nieder80} gave a formula of evaluating  $|\L_{k,r}(0,0;m,n)|$ for positive integers $k,r$. In fact, Niederhausen considered a multi-dimensional lattice path enumeration problem, and his method is based on the finite operator calculus developed in Rota et al. \cite{RKO}. It should be mentioned that B\"{o}hm \cite[(10)]{Bohm2010} studied an enumeration problem which is equivalent to the evaluation of $|\L_{k,r}(0,0;m,n)|$.
This equivalence will be illustrated in Section \ref{sec-re}. A formula of $|\L_{\frac{1}{k},\frac{c}{k}}(0,0;m,n)|$ with positive integers $k,c$
was given in Niederhausen \cite[Proposition 4]{Nieder02}.

Among the above lattice path enumeration problems, the starting points are fixed to be the origin of the plane. For this reason, all manipulations on the lattice paths are naturally required to fix the starting points. In fact, all known methods obeyed this spirit. Since we are considering the enumeration of lattice paths starting from a point not necessarily being fixed (usually the origin is chosen as the starting point), we could use induction on the starting points. This suggests a forward recursion approach to the enumeration problem considered here.

This paper is organized as follows.
In Section \ref{sec-2}, by using the forward recursion, we obtain a formula of $|L_{k,r}(a,b;m,n)|$ for integers $k,r,a,b,m,n$, as well as that of  $|\L_{k,r}(a,b;m,n)|$, $|L_{\frac{1}{k},r}(a,b;m,n)|$ and $|\L_{\frac{1}{k},r}(a,b;m,n)|$.
In Section \ref{sec-re}, we shall show how to derive
Koroljuk's formula from the results in Section \ref{sec-2}.
In particular, we shall show the relations among the lattice paths with unit steps $(1,1)$ and $(-p,1)$ studied by Koroljuk, the lattice paths with unit steps $(1,t)$ and $(1,-1)$ studied by B\"{o}hm, and the lattice paths with unit steps $(1,0)$ and $(0,1)$ studied by Niederhausen, where $p,t$ are some fixed positive integers.

\section{The Lattice Path Enumeration Formulas}\label{sec-2}

In this section we give an enumeration formula for $|L_{k,r}(a,b;m,n)|$. By using some simple bijections, we then derive the enumeration formulas for $|\L_{k,r}(a,b;m,n)|$, $|L_{\frac{1}{k},r}(a,b;m,n)|$ and $|\L_{\frac{1}{k},r}(a,b;m,n)|$ from that for $|L_{k,r}(a,b;m,n)|$.

Our first main result is as follows.

\begin{theo}\label{lat-abmnkc} Suppose that $k, r, a, b, m$ and
$n$ are integers satisfying
\begin{align}\label{eq-co2.1}
k\geq 1,\quad  0\leq a\leq m, \quad n\geq km-r, \quad \max\{0,ka-r\}\leq b\leq n.
\end{align}
Then we have
\begin{eqnarray}\label{eq-Lkrabmn}
&|L_{k,r}(a,b;m,n)|=\sum\limits_{i=0}^{\lfloor\frac{b+r-ka}{k+1}\rfloor}(-1)^i
  \frac{n+r+1-km}{n+r+1-k(a+i)}{m+n+r-(k+1)(a+i)\choose m-a-i}{b+r-k(a+i)\choose i}.
\end{eqnarray}
\end{theo}

\noindent{\bf Remark.}
Note that when $a=b=r=0$, $k\rightarrow k-1$ and $n=(k-1)m$, \eqref{eq-Lkrabmn} becomes $$|L_{k-1,0}(0,0;m,(k-1)m)|=\frac{1}{(k-1)m+1}{km\choose m},$$
which are known as the order-$k$ Fuss-Catalan numbers, see \cite{BacKra11, Fuss91} and \cite[Eq. (7.67)]{G-K-P}. They also enumerate $(k-1)$-Dyck paths, as well as $k$-ary trees. In particular, when $k=2$ they are the well known Catalan numbers, and when $k=3$ they count $2$-Dyck paths, as well as ternary trees. For more details, see \cite{Aval08, Stanley1999}.

We shall prove Theorem \ref{lat-abmnkc} by induction. Before doing this, let us first give two lemmas, one of which is concerned with a recursion relation on $|L_{k,r}(a,b;m,n)|$, and the other some initial values of  $|L_{k,r}(a,b;m,n)|$. For the recursion relation on $|L_{k,r}(a,b;m,n)|$, we have the following result.

\begin{lem}\label{prop-Nabmna-u}
Suppose that $k, r, a, b, m$ and
$n$ are integers satisfying
\begin{align}\label{cond-1}
k,m\geq 1,\quad  0\leq a\leq m, \quad n\geq km-r,  \quad k(a+1)-r\leq b\leq n.
\end{align}
Then we have
\begin{align*}
 |L_{k,r}(a,b+1;m,n)|=|L_{k,r}(a,b;m,n)|-|L_{k,r}(a,b-k;m-1,n-k)|.
\end{align*}
\end{lem}

\pf Note that, for each lattice path of $L_{k,r}(a,b;m,n)$, the first step is either from $(a,b)$ to $(a+1,b)$ or from $(a,b)$ to $(a,b+1)$.
Moreover, if $k, r, a, b, m$ and
$n$ satisfy \eqref{cond-1}, then $(a,b)$ and $(m,n)$ must stay above the boundary line $y=kx-r$. It is easy to see that
\begin{align*}
|L_{k,r}(a,b;m,n)|=|L_{k,r}(a+1,b;m,n)|+|L_{k,r}(a,b+1;m,n)|.
\end{align*}
Now it suffices to show that
\begin{equation}\label{eq-2.3}
|L_{k,r}(a+1,b;m,n)|=|L_{k,r}(a,b-k;m-1,n-k)|.
\end{equation}
To this end, define a map
$$
f:L_{k,r}(a+1,b;m,n)\longrightarrow L_{k,r}(a,b-k;m-1,n-k)
$$
as follows: for each lattice path $\xi \in L_{k,r}(a+1,b;m,n)$, let $f(\xi)$ be the lattice path obtained from $\xi$ by translating $\xi$  horizontally left by one unit and vertically down by $k$ units.
Thus, the resulting path $f(\xi)$ is from $(a,b-k)$ to $(m-1,n-k)$.
We claim that $f(\xi)\in L_{k,r}(a,b-k;m-1,n-k)$, since each lattice point $(x,y)\in \xi$ satisfies $y\geq kx-r$, and its image $(x-1,y-k)\in f(\xi)$ naturally satisfies $y-k\geq  k(x-1)-r$, i.e., $f(\xi)$ stays above the line $y=kx-r$. Clearly, the map $f$ is reversible and unique. Therefore, we obtain the equality \eqref{eq-2.3}. This completes the proof.
\qed

Next we give a formula for $|L_{k,r}(a,b;m,n)|$ when $r=0$ and $0\leq b-ka\leq k$. For convenience, we abbreviate $L_{k,0}(a,b;m,n)$ by $L_k(a,b;m,n)$. The following result serves as the base case of the induction proof of Theorem \ref{lat-abmnkc}.

\begin{lem}\label{Nabmna-ua-1}
Suppose that $k, a, b, m$ and
$n$ are integers satisfying
$$k,m\geq 1,\quad  0\leq a\leq m, \quad 0\leq b\leq n, \quad n\geq km, \quad 0\leq b-ka\leq k.$$
Then we have
\begin{align*}
|L_k(a,b;m,n)|=\frac{n+1-km}{n+1-ka}{m+n-(k+1)a\choose m-a}.
\end{align*}
\end{lem}

\pf
Given the conditions satisfied by $a, b, m, n$ and $k$, each path in $L_k(a,b;m,n)$ must pass through the point $(a,ka+k)$. Thus, we have
$$
|L_k(a,b;m,n)|=|L_k(a,ka;m,n)|,
$$
since the path from $(a,b)$ to $(a,ka+k)$ is unique due to the restriction of unit steps.

Now we are to compute $|L_k(a,ka;m,n)|$. Let $\C_k(a,b;m,n)$ denote the set of lattice paths from $(a,b)$ to $(m,n)$ which cross the boundary line $y=kx$. Since the number of lattice paths from $(a,ka)$ to $(m,n)$ with unit steps $(1,0)$ and $(0,1)$ is $\binom{m+n-(k+1)a}{m-a}$, we get that
\begin{equation*}
|L_k(a,ka;m,n)|=\binom{m+n-(k+1)a}{m-a}-|\C_k(a,ka;m,n)|.
\end{equation*}
By translating each lattice path in $\C_k(a,ka;m,n)$ horizontally left by $a$ units and vertically down by $ka$ units, we shall obtain a lattice path in $\C_k(0,0;m-a,n-ka)$. Conversely, each path in $\C_k(0,0;m-a,n-ka)$ can be obtained in such a way. Therefore, we have
$$|\C_k(a,ka;m,n)|=|\C_k(0,0;m-a,n-ka)|.$$
By the solution of the generalized ballot problem \cite{Gou-Ser}, we see that
\begin{align*}
|\C_k(0,0;m,n)|=k{m+n\choose m-1},
\end{align*}
and hence
\begin{align*}
|\C_k(0,0;m-a,n-ka)|=k{m+n-(k+1)a\choose m-a-1}.
\end{align*}
Therefore,
\begin{align*}
|L_k(a,b;m,n)|&=\binom{m+n-(k+1)a}{m-a}-k{m+n-(k+1)a\choose m-a-1}\\[5pt]
&=\frac{n+1-km}{n+1-ka}{m+n-(k+1)a\choose m-a}.
\end{align*}
This completes the proof. \qed

We now proceed to prove Theorem \ref{lat-abmnkc}.

\noindent {\it Proof of Theorem \ref{lat-abmnkc}.}
There is a one-to-one correspondence between $L_{k,r}(a,b;m,n)$ and $L_k(a,b+r;m,n+r)$, since a lattice path in $L_{k,r}(a,b;m,n)$ will be transformed into a lattice path in $L_k(a,b+r;m,n+r)$ by simply shifting the $x$-axis vertically down or up by $|r|$ units according to the sign of $r$. This means that
\begin{equation}\label{eq-relr0}
|L_{k,r}(a,b;m,n)|=|L_k(a,b+r;m,n+r)|.
\end{equation}
Now it is sufficient to prove the theorem for the case of $r=0$, namely
\begin{equation}\label{enu-labmna}
|L_k(a,b;m,n)|=\sum\limits_{i=0}^{\lfloor\frac{b-ka}{k+1}\rfloor}(-1)^i
  \frac{n+1-km}{n+1-k(a+i)}{m+n-(k+1)(a+i)\choose m-a-i}{b-k(a+i)\choose i}.
\end{equation}
For $r\neq 0$, one can easily derive \eqref{eq-Lkrabmn} from \eqref{enu-labmna} together with \eqref{eq-relr0}.

In the following we shall prove \eqref{enu-labmna} by induction on the value $\ell(a,b,k)=\lfloor\frac{b-ka}{k+1}\rfloor$. For $\ell(a,b,k)=0$, namely, $0\leq b-ka\leq k$, it is easy to see that \eqref{enu-labmna} holds true by Lemma \ref{Nabmna-ua-1}. Assume that \eqref{enu-labmna} holds for any $a,b,k$ satisfying $\lfloor\frac{b-ka}{k+1}\rfloor<\ell$. We are to prove that it is also valid for $b=\ell(k+1)+ka+j$ for some $0\leq j\leq k$, namely
\begin{equation}\label{eq-6.21e}
|L_k(a,\ell(k+1)+ka+j;m,n)|
=\sum\limits_{i=0}^{\ell}(-1)^i {\ell(k+1)+j-ki\choose i} h(a,m,n,k,i),
\end{equation}
where
\begin{equation*}
h(a,m,n,k,i)=\frac{n+1-km}{n+1-k(a+i)}{m+n-(k+1)(a+i)\choose m-a-i}.
\end{equation*}
By Lemma \ref{prop-Nabmna-u}, we have
\begin{align*}
|L_k(a,b;m,n)|=|L_k(a,b-1;m,n)|-|L_k(a,b-1-k;m-1,n-k)|.
\end{align*}
Since for $L_k(a,b-1-k;m-1,n-k)$ it holds
$$\ell(a,b-1-k,k)=\left\lfloor\frac{(b-1-k)-ka}{k+1}\right\rfloor=
\left\lfloor\frac{\ell(k+1)+ka+j-1-k-ka}{k+1}\right\rfloor=\ell-1,$$
we can determine the value of $|L_k(a,b-1-k;m-1,n-k)|$ by the induction hypothesis. However, for $L_k(a,b-1;m,n)$ it holds
$$\ell(a,b-1,k)=\left\lfloor\frac{(b-1)-ka}{k+1}\right\rfloor
=\left\lfloor\frac{\ell(k+1)+ka+j-1-ka}{k+1}\right\rfloor=\ell$$
for $j\geq 1$. In order to use the induction hypothesis, we need to apply Lemma \ref{prop-Nabmna-u} to $L_k(a,b-1;m,n)$, and we get
\begin{align*}
|L_k(a,b-1;m,n)|=|L_k(a,b-2;m,n)|-|L_k(a,b-2-k;m-1,n-k)|.
\end{align*}
Now we have
\begin{align*}
|L_k(a,b;m,n)|&=|L_k(a,b-2;m,n)|-\sum\limits_{t=0}^{1}
|L_k(a,b-2+t-k;m-1,n-k)|.
\end{align*}
Again, we check that $L_k(a,b-2-k;m-1,n-k)$ satisfies the induction hypothesis, while $L_k(a,b-2;m,n)$ does not for $j\geq 2$. We now need to
apply Lemma \ref{prop-Nabmna-u} to $L_k(a,b-2;m,n)$. Iterate this procedure until reaching the following formula
\begin{align*}
|L_k(a,b;m,n)|&=|L_k(a,b-j-1;m,n)|-\sum\limits_{t=0}^{j}|L_k(a,b-j-1+t-k;m-1,n-k)|.
\end{align*}
By using the induction hypothesis and substituting $\ell(k+1)+ka+j$ for $b$, we find that
\begin{align*}
|L_k(a,b-j-1;m,n)|
=\sum\limits_{i=0}^{\ell-1}(-1)^i {\ell(k+1)-1-ki\choose i}h(a,m,n,k,i),
\end{align*}
and
\begin{align*}
& \sum\limits_{t=0}^j|L_k(a,b-j-1+t-k;m-1,n-k)|\nonumber\\
=&\sum\limits_{t=0}^j\sum\limits_{i=0}^{\ell-1}(-1)^i {(\ell-1)(k+1)+t-ki\choose i}h(a,m-1,n-k,k,i)\nonumber\\
=&\sum\limits_{t=0}^j\sum\limits_{i=1}^{\ell}(-1)^{i-1}{\ell(k+1)+t-1-ki\choose i-1}h(a,m,n,k,i),
\end{align*}
where the last equality follows from the fact that $h(a,m-1,n-k,k,i-1)=h(a,m,n,k,i)$.
Therefore,
\begin{align*}
|L_k(a,b;m,n)|=&\sum\limits_{i=0}^{\ell-1}(-1)^i {\ell(k+1)-1-ki\choose i}h(a,m,n,k,i)\\[6pt]
&-\sum\limits_{t=0}^j\sum\limits_{i=1}^{\ell}(-1)^{i-1}{\ell(k+1)+t-1-ki\choose i-1}h(a,m,n,k,i)\\[6pt]
=&\sum\limits_{i=0}^{\ell-1}(-1)^i {\ell(k+1)-1-ki\choose i}h(a,m,n,k,i)\\
&-\sum\limits_{i=1}^{\ell}(-1)^{i-1}\left(\sum\limits_{t=0}^j{\ell(k+1)+t-1-ki\choose i-1}\right)h(a,m,n,k,i).
\end{align*}
Note that
\begin{align*}
\sum\limits_{t=0}^j{\ell(k+1)+t-1-ki\choose i-1}&=\sum\limits_{t=0}^j\left({\ell(k+1)+t-ki\choose i}-{\ell(k+1)+t-1-ki\choose i}\right)\\[6pt]
&={\ell(k+1)+j-ki\choose i}-{\ell(k+1)-1-ki\choose i}.
\end{align*}

Thus,
\begin{align*}
|L_k(a,b;m,n)|
=&\sum\limits_{i=0}^{\ell-1}(-1)^i {\ell(k+1)-1-ki\choose i}h(a,m,n,k,i)\\[6pt]
&-\sum\limits_{i=1}^{\ell}(-1)^{i}{\ell(k+1)-1-ki\choose i}h(a,m,n,k,i)\\[6pt]
&-\sum\limits_{i=1}^{\ell}(-1)^{i-1}{\ell(k+1)+j-ki\choose i}h(a,m,n,k,i)\\[6pt]
=&h(a,m,n,k,0)+\sum\limits_{i=1}^{\ell}(-1)^{i}{\ell(k+1)+j-ki\choose i}h(a,m,n,k,i)\\[6pt]
=&\sum\limits_{i=0}^{\ell}(-1)^{i}{\ell(k+1)+j-ki\choose i}h(a,m,n,k,i),
\end{align*}
as desired. This completes the proof.
\qed

Now we are able to give a formula to compute $|\L_{k,r}(a,b;m,n)|$.

\begin{coro}\label{coro-2.1}
Suppose that $k, r, a, b, m$ and
$n$ are integers satisfying
$$k\geq 1,\quad  0\leq a\leq m, \quad \max\{0,ka-r\}<b\leq n, \quad n>km-r.$$
Then we have
{\small
\begin{eqnarray}\label{eq-snL}
&|\L_{k,r}(a,b;m,n)|=
\sum\limits_{i=0}^{\lfloor\frac{b+r-1-ka}{k+1}\rfloor}(-1)^i
  \frac{n+r-km}{m+n+r-(k+1)(a+i)}{m+n+r-(k+1)(a+i)\choose m-a-i} {b+r-1-k(a+i)\choose i}.
\end{eqnarray}
}
\end{coro}

\pf
There is an obvious bijection between $\L_{k,r}(a,b;m,n)$ and $L_{k,r}(a,b-1;m,n-1)$. Shift any path $\xi\in\L_{k,r}(a,b;m,n)$ by one vertical down unit, and the resulting path  obviously belongs to $L_{k,r}(a,b-1;m,n-1)$. Clearly, this map is reversible and unique. Thus
\begin{align*}
|\L_{k,r}(a,b;m,n)|=|L_{k,r}(a,b-1;m,n-1)|.
\end{align*}
By \eqref{eq-Lkrabmn}, we immediately obtain the desired result. This completes the proof.
\qed

Theorem \ref{lat-abmnkc} can also be used to derive the number of lattice paths above the boundary line $y=\frac{x}{k}-r$ for some positive integer $k$.

\begin{coro}\label{coro-2.2}
Suppose that $k, r, a, b, m$ and $n$ are integers satisfying
$$k\geq 1,\quad  0\leq a\leq m, \quad \max\{0,\frac{a}{k}-r\}\leq b\leq n, \quad n\geq\frac{m}{k}-r.$$
Then we have
\begin{eqnarray}\label{eq-3.1}
&|L_{\frac{1}{k},r}(a,b;m,n)|=
\sum\limits_{i=0}^{\lfloor\frac{k(n+r)-m}{k+1}\rfloor}(-1)^i
  \frac{k(b+r)-a+1}{k(n+r-i)-a+1}{(k+1)(n-i)-a-b+kr\choose n-b-i} {k(n+r-i)-m\choose i}.
\end{eqnarray}
\end{coro}
\pf
We shall show that there is a one-to-one correspondence between $L_{\frac{1}{k},r}(a,b;m,n)$ and $L_k(0,k(n+r)-m;n-b,k(n+r)-a)$.
Define a map $f: L_{\frac{1}{k},r}(a,b;m,n)\rightarrow L_k(0,k(n+r)-m;n-b,k(n+r)-a)$ in the following two steps.
\begin{itemize}
\item[1.] First, shift the $x$-axis vertically up by $n$ units and the $y$-axis horizontally right by $k(n+r)$ units such that the new origin is located at $(k(n+r),n)$, and then reverse the direction of each axis. Now, a lattice path $\xi\in L_{\frac{1}{k},r}(a,b;m,n)$ in the original coordinate system becomes a lattice path from $(k(n+r)-m,0)$ to $(k(n+r)-a,n-b)$ that stays below the line $y=\frac{x}{k}$ in the new coordinate system, denoted by $\xi_1$.

\item[2.] Second, reflect $\xi_1$ along the diagonal $y=x$. Denote the resulting path by $f(\xi)$. It is clear that $f(\xi)$ is a lattice path from $(0,k(n+r)-m)$ to $(n-b,k(n+r)-a)$ that stays above the line $y=kx$, namely, $f(\xi)\in L_k(0,k(n+r)-m;n-b,k(n+r)-a)$.
\end{itemize}
Clearly, each step is reversible and unique. Therefore, $f$ is bijective.
We immediately derive the desired result from \eqref{enu-labmna}. This completes the proof.
\qed

With the aid of  Corollary \ref{coro-2.2}, we
obtain the following result for counting the lattice paths of $\L_{\frac{1}{k},r}(a,b;m,n)$. The proof is similar to
that of Corollary \ref{coro-2.2} as derived from Theorem \ref{lat-abmnkc} and is omitted.

\begin{coro}\label{coro-2.3}
Suppose that $k, r, a, b, m$ and $n$ are integers satisfying
$$k\geq 1,\quad  0\leq a\leq m, \quad \max\{0,\frac{a}{k}-r\}<b\leq n, \quad n>\frac{m}{k}-r.$$
Then we have
{\small
\begin{eqnarray}\label{eq-3.3}
&\hskip -9mm|\L_{\frac{1}{k},r}(a,b;m,n)|
=\sum\limits_{i=0}^{\lfloor\frac{k(n+r)-m-1}{k+1}\rfloor}(-1)^i
   \frac{k(b+r)-a}{(k+1)(n-i)-a-b+kr}{(k+1)(n-i)-a-b+kr\choose n-b-i}{k(n+r-i)-m-1\choose i}.
\end{eqnarray}
}
\end{coro}

It should be mentioned that Niederhausen \cite[Proposition 4]{Nieder02} gave a formula for the number of lattice paths from $(0,0)$ to $(n,m)$ staying strictly above the line $y=(x-c)/k$ for positive integers $k$ and $c$, by using Koroljuk's method of decomposing lattice paths.
This enumeration formula is equivalent to $|\L_{\frac{1}{k},\frac{c}{k}}(0,0;n,m)|$.

\noindent{\bf Remark.}
Note that when $r$ is a real number but not an integer, the restriction of the boundary $y=kx-r$ on the lattice paths is equivalent to that of the line $y=kx-\lfloor r\rfloor$ for \eqref{eq-Lkrabmn}. Now we have $|L_{k,r}(a,b;m,n)|=|\L_{k,r}(a,b;m,n)|=|L_{k,\lfloor r\rfloor}(a,b;m,n)|$.
In addition, as long as $kr$ is an integer, the formulae \eqref{eq-3.1} and \eqref{eq-3.3} still work. If $kr$ is not an integer, it is easy to see that $L_{\frac{1}{k},r}(a,b;m,n)=\L_{\frac{1}{k},r}(a,b;m,n)=L_{\frac{1}{k},\frac{\lfloor kr\rfloor}{k}}(a,b;m,n)$. Substituting $\lfloor kr\rfloor$ for $kr$ in \eqref{eq-3.1}, the formula still works.

\section{Koroljuk's formula}\label{sec-re}

The aim of this section is to give a new proof of Koroljuk's formula \eqref{eq-Koro} by using the above enumeration results. To this end, we first prove the following result, a counterpart of Koroljuk's formula.

\begin{theo}\label{theo-3.1}
For some positive integers $m,n,p$ and $c$,  the number of lattice paths from $(0,0)$ to $(m-pn,m+n)$, using the steps $(1,1)$ and $(-p,1)$ in the plane and not intersecting the line $x=c$, is enumerated by $|\L_{p,c+pn-m}(0,0;n,m)|$.
\end{theo}

\pf
Given a lattice path $\xi$ from $(0,0)$ to $(m-pn,m+n)$, using the steps $(1,1)$ and $(-p,1)$ in the plane and not intersecting the line $x=c$, we shall map it into a lattice path in $\L_{p,c+pn-m}(0,0;n,m)$ in the following way:
\begin{itemize}
\item[1.] First, rotate $\xi$ by $90^\circ$ clockwise about the origin and then shift the $x$-axis by $c$ units downwards. Denote the image by $\xi_0$ (see Figure \ref{fig-a}), which is a path from $(0,c)$ to $(m+n,c+pn-m)$ with $(1,p),(1,-1)$-steps that stays strictly above the $x$-axis. Let $v=c+pn-m$, which is obviously a positive integer.

\item[2.] Second, rotate $\xi_0$ by $180^\circ$ about the origin, and then shift the image horizontally right such that the leftmost lattice point is located on the $y$-axis. The resulting image, denoted by $\xi_1$  (see Figure \ref{fig-b}), is a path from $(0,-v)$ to $(m+n,-c)$ with $(1,p),(1,-1)$-steps that stays strictly below the line $y=0$.

\item[3.] Thirdly, shift the $y$-axis horizontally left by $v$ units. Then rotate the axes by $45^\circ$ counterclockwise about the origin. Next reverse the direction of the $y$-axis and then exchange the two axes. Now scale the Cartesian coordinate system with a scale factor $\cos\frac{\pi}{4}$. The resulting path, denoted by $\xi_2$  (see Figure \ref{fig-c}), is from $(v,0)$ to $(c+\frac{\sqrt{2+2r^2}}{2}n\sin\theta ,\frac{\sqrt{2+2r^2}}{2}n\sin\theta)$, where $\theta=\frac{3\pi}{4}-\arctan p$ is the angle between the non-horizontal step and the horizontal step. Moreover, $\xi_2$ stays strictly below the line $y=x$.

\item[4.] Next, keep the portion of $\xi_2$ on the $x$-axis fixed and shear the remaining portion of $\xi_2$ rightward parallel to the $x$-axis such that all the non-horizontal steps become vertical, which means there is a rotation by $\beta_1=\frac{\pi}{2}-\theta$ degrees clockwise for each non-horizontal step. Then scale the length of each vertical step to $1$. At the same time, rotate the boundary line $y=x$ by $\beta_2=\arctan p-\frac{\pi}{4}$ degrees clockwise about the origin. The boundary line becomes $y=x/p$. Denote the resulting path by $\xi_3$  (see Figure \ref{fig-d}). Note that $\beta_1=\beta_2$. It is easy to see that $\xi_3$ goes from $(v,0)$ to $(c+pn,n)$ and stays strictly below the line $y=x/p$.

\item[5.] Finally, shift the $y$-axis horizontally right by $v$ units, and then reflect $\xi_3$ along the diagonal $y=x$. Denote the image by $\xi'$. It is clear that $\xi'\in\L_{p,c+pn-m}(0,0;n,m)$.
\end{itemize}
\vskip -7pt
\begin{figure}[!ht]
\setlength\unitlength{2.5pt}
\begin{minipage}[!h]{0.49\textwidth}
\centering
\begin{picture}(74,42)
\put(2,8){\vector(1,0){68}}
\put(68,5.2){\small $x$}
\put(8,2){\vector(0,1){35}}
\put(5.27,35.5){\small $y$}
\put(5.6,5.0){\footnotesize $0$}
\multiput(8,18)(24,0){2}{\circle*{1.2}}
\multiput(12,26)(4,-4){5}{\circle*{1.2}}
\multiput(36,26)(4,-4){4}{\circle*{1.2}}
\multiput(52,22)(4,-4){3}{\circle*{1.2}}
\put(5.1,17){\small $c$}
\thicklines
\put(8,18){\line(1,2){4}}
\put(12,26){\line(1,-1){16}}
\put(28,10){\line(1,2){8}}
\put(36,26){\line(1,-1){12}}
\put(48,14){\line(1,2){4}}
\put(52,22){\line(1,-1){8}}
\put(63,9.7){\small $v$}
\put(62,11){\makebox(0,0){$\left\} \rule[0mm]{0mm}{3.8mm} \right.$}}
\end{picture}
\vskip -3pt
\caption{$\xi_0$ with $p=2$ and $n=4$. \label{fig-a}}
\end{minipage}
\hspace{0.01\linewidth}
\setlength\unitlength{2.5pt}
\begin{minipage}[!h]{0.49\textwidth}
\centering
\begin{picture}(78,40)
\put(4.5,26){\vector(1,0){68}}
\put(70.2,23.3){\small $x$}
\put(10,3){\vector(0,1){33.6}}
\put(7.27,35.1){\small $y$}
\put(8.1,23.3){\scriptsize $0$}
\put(64.1,21){\makebox(0,0){$\left\} \rule[0mm]{0mm}{5.26mm} \right.$}}
\put(3.6,19.3){\small $-v$}
\put(65.3,20.3){\small $c$}
\multiput(10,20)(4,-4){3}{\circle*{1.2}}
\multiput(22,20)(4,-4){4}{\circle*{1.2}}
\multiput(42,24)(4,-4){5}{\circle*{1.2}}
\multiput(38,16)(24,0){2}{\circle*{1.2}}
\thicklines
\put(10,20){\line(1,-1){8}}
\put(18,12){\line(1,2){4}}
\put(22,20){\line(1,-1){12}}
\put(34,8){\line(1,2){8}}
\put(42,24){\line(1,-1){16}}
\put(58,8){\line(1,2){4}}
\end{picture}
\caption{$\xi_1$ with $p=2$ and $n=4$.  \label{fig-b}}
\end{minipage}\\
\vskip -5pt
\begin{minipage}[!h]{0.5\textwidth}
\setlength\unitlength{3.1pt}
\centering
\begin{picture}(68,50)
\thinlines
\put(3,8){\vector(1,0){58}}
\put(8,3){\vector(0,1){38}}
\put(3.9,3.9){\line(1,1){36.6}}
\put(59.2,5.5){\small $x$}
\put(5.5,39.5){\small $y$}
\put(8.15,5.2){\footnotesize $0$}
\put(14.65,4.9){\small $v$}
\multiput(32,32)(2,0){6}{\line(1,0){1}}
\put(38,30.5){\makebox(0,0){$\underbrace{\rule[0mm]{13.5mm}{0mm}}$}}
\put(37.5,27.6){\footnotesize $c$}
\put(28,37){\small $y=x$}
\multiput(16,8)(4,0){3}{\circle*{1.2}}
\multiput(22,14)(4,0){4}{\circle*{1.2}}
\multiput(30,26)(4,0){5}{\circle*{1.2}}
\multiput(32,20)(12,12){2}{\circle*{1.2}}
\qbezier(22.5,8)(22.7,9)(23.5,9.2)
\put(21.5,8.76){\scriptsize $\theta$}
\qbezier(32.5,14)(32.7,15)(33.5,15.2)
\put(31.2,15.15){\scriptsize $\theta$}
\qbezier(44.5,26)(44.7,27)(45.5,27.2)
\put(43.2,27.15){\scriptsize $\theta$}
\thicklines
\put(16,8){\line(1,0){8}}
\put(24,8){\line(-1,3){2}}
\put(22,14){\line(1,0){12}}
\put(34,14){\line(-1,3){4}}
\put(30,26){\line(1,0){16}}
\put(46,26){\line(-1,3){2}}
\end{picture}
\vskip -3mm
\caption{$\xi_2$ with $p=2$ and $n=4$. \label{fig-c}}
\end{minipage}
\hspace{0.01\linewidth}
\begin{minipage}[!h]{0.49\textwidth}
\centering
\setlength\unitlength{3.1pt}
\begin{picture}(70,51)
\put(3,8){\vector(1,0){59}}
\put(8,3){\vector(0,1){38}}
\put(4,6){\line(2,1){55}}
\put(60.2,5.5){\small $x$}
\put(5.5,39.5){\small $y$}
\put(8.15,5.1){\footnotesize $0$}
\put(14.65,4.9){\small $v$}
\put(52,8){\line(0,1){0.75}}
\put(8,24){\line(1,0){0.75}}
\put(47.6,5.2){\small $m$$+$$v$}
\put(5.7,23.5){\small $n$}
\put(45.0,33.0){\small $y=\frac{x}{p}$}
\multiput(16,8)(4,0){3}{\circle*{1.2}}
\multiput(24,12)(4,0){4}{\circle*{1.2}}
\multiput(36,20)(4,0){5}{\circle*{1.2}}
\multiput(36,16)(16,8){2}{\circle*{1.2}}
\thicklines
\put(16,8){\line(1,0){8}}
\put(24,8){\line(0,1){4}}
\put(24,12){\line(1,0){12}}
\put(36,12){\line(0,1){8}}
\put(36,20){\line(1,0){16}}
\put(52,20){\line(0,1){4}}
\end{picture}
\vskip -3mm
\caption{$\xi_3$  with $p=2$ and $n=4$. \label{fig-d}}
\end{minipage}
\end{figure}
\vskip 3pt

It is easy to verify that each step is reversible and unique. Therefore the map is a bijection. This completes the proof. \qed

\noindent \textbf{Remark.} The lattice paths obtained in the first step of the proof have been studied by B\"{o}hm \cite{Bohm2010}. As pointed out by B\"{o}hm \cite[Eq. (10)]{Bohm2010}, the number of paths, going from initial altitude $m>0$ to terminal altitude $k>0$ with $(1,r),(1,-1)$-steps such that the paths do not touch or cross the $x$-axis and the number of $(1,r)$-steps is equal to $a$, can be given by
\begin{equation}\label{eq-NieBoh}
N_{m,k}(a)=\sum\limits_{\ell\geq 0}(-1)^{\ell}\frac{m}{m+(r+1)(a-\ell)}{m+(r+1)(a-\ell)\choose a-\ell}{k-r\ell-1\choose \ell}.
\end{equation}
B\"{o}hm \cite{Bohm2010} attributes the above formula to Niederhausen. By the proof of the theorem, we see that $N_{m,k}(a)$ also counts the lattice paths from $(0,0)$ to $(a,m+ra-k)$ that strictly stay above the line $y=rx-k$. Thus, $N_{m,k}(a)=|\L_{r,k}(0,0;a,m+ra-k)|$. By \eqref{eq-snL}, we can also obtain Niederhausen's formula \eqref{eq-NieBoh}.

Now we can prove Koroljuk's formula.

\begin{coro}
For some positive integers $m,n,p$ and $c$,  the number of lattice paths from $(0,0)$ to $(m-pn,m+n)$, using the steps $(1,1)$ and $(-p,1)$ in the plane and intersecting the line $x=c$, is enumerated by \eqref{eq-Koro}.
\end{coro}

\pf Note that there are $\binom{m+n}{n}$ lattice paths from $(0,0)$ to $(m-pn,m+n)$ by using the steps $(1,1)$ and $(-p,1)$. By the above theorem, the number of lattice paths from $(0,0)$ to $(m-pn,m+n)$, using the steps $(1,1)$ and $(-p,1)$ in the plane and intersecting the line $x=c$, is counted by
\begin{align*}
\binom{m+n}{n}-|\L_{p,c+pn-m}(0,0;n,m)|.
\end{align*}
By \eqref{eq-snL}, it is equal to
\begin{eqnarray}\label{eq-3.7}
&\binom{m+n}{n}-
\sum\limits_{i=0}^{\lfloor\frac{c+pn-m-1}{p+1}\rfloor}(-1)^i
  \frac{c}{c+(p+1)(n-i)}{c+(p+1)(n-i)\choose n-i} {c-m-1+p(n-i)\choose i}.
\end{eqnarray}
It is well known that
\begin{align*}
{x\choose k}=(-1)^k{k-x-1\choose k},\quad\quad \mbox{ for } x\in\R,\  k\in\mathbb{N},
\end{align*}
see \cite[Eq. (17), p. 58]{Knuth}. Thus, we have
$$
{c-m-1+p(n-i)\choose i}=(-1)^i {m+n-c-(p+1)(n-i)\choose i},
$$
and \eqref{eq-3.7} is equal to
\begin{eqnarray*}
&\binom{m+n}{n}-
\sum\limits_{i=0}^{\lfloor\frac{c+pn-m-1}{p+1}\rfloor}
  \frac{c}{c+(p+1)(n-i)}{c+(p+1)(n-i)\choose n-i} {m+n-c-(p+1)(n-i)\choose i}\\[5pt]
=&\binom{m+n}{n}-
\sum\limits_{i=n-\lfloor\frac{c+pn-m-1}{p+1}\rfloor}^{n}
  \frac{c}{c+(p+1)i}{c+(p+1)i\choose i} {m+n-c-(p+1)i\choose n-i}.
\end{eqnarray*}
By the following Hagen-Rothe identity \cite[Eq. (11)]{Gould},
\begin{align*}
\sum\limits_{i=0}^n\frac{\gamma}{\gamma+\beta i}{\gamma+\beta i\choose i}{\alpha+\beta(n-i)\choose n-i}={\alpha+\gamma+\beta n\choose n},
\quad \alpha,\beta,\gamma\in\R,\ n\in\mathbb{N},
\end{align*}
we find that
$$
\sum\limits_{i=0}^{n}
  \frac{c}{c+(p+1)i}{c+(p+1)i\choose i} {m+n-c-(p+1)i\choose n-i}=\binom{m+n}{n}.
$$
Therefore, \eqref{eq-3.7} is equal to
$$
\sum\limits_{i=0}^{n-\lfloor\frac{c+pn-m-1}{p+1}\rfloor-1}
  \frac{c}{c+(p+1)i}{c+(p+1)i\choose i} {m+n-c-(p+1)i\choose n-i},
$$
that is
\begin{align}\label{ref-koro}
\sum\limits_{i=0}^{\lfloor\frac{m+n-c}{p+1}\rfloor}\frac{c}{c+(p+1)i}{c+(p+1)i\choose i}{m+n-c-(p+1)i\choose n-i},
\end{align}
since, for positive integers $m,n,p$ and $c$, we have
\begin{align*}
n-\left\lfloor\frac{c+pn-m-1}{p+1}\right\rfloor-1&=n+\left\lceil\frac{m-pn-c+1}{p+1}\right\rceil-1\\[5pt]
&=\left\lceil \frac{m+n-c+1}{p+1}\right\rceil-1\\[5pt]
&=\left\lfloor\frac{m+n-c}{p+1}\right\rfloor.
\end{align*}
By the substitution of $i=\frac{s-c}{p+1}$ for the summation index, one can obtain Koroljuk's formula \eqref{eq-Koro} from \eqref{ref-koro} immediately. This completes the proof.
\qed

As noted in the introduction, Niederhausen \cite[p. 9]{Nieder02} studied the enumeration of lattice paths from $(0,0)$ to $(m,n)$ with $(1,0), (0,1)$-steps that stay strictly  above the line $y=k(x-d)$, an equivalent problem to Koroljuk's formula. By \eqref{eq-snL}, such lattice paths are counted by
\begin{align*}
|\L_{k,kd}&(0,0;m,n)|\\
=&\sum\limits_{i=0}^{\lfloor\frac{kd-1}{k+1}\rfloor}(-1)^i
  \frac{n+kd-km}{m+n+kd-(k+1)i}{m+n+kd-(k+1)i\choose m-i} {kd-1-ki\choose i}\\
=&\sum\limits_{i=0}^{\lfloor\frac{kd-1}{k+1}\rfloor}
  \frac{n+kd-km}{m+n+kd-(k+1)i}{m+n+kd-(k+1)i\choose m-i} {ki-kd+i\choose i}\\
=&\sum\limits_{i=0}^{m}
  \frac{n+kd-km}{m+n+kd-(k+1)i}{m+n+kd-(k+1)i\choose m-i} {ki-kd+i\choose i}\\
&-\sum\limits_{i=\lfloor\frac{kd-1}{k+1}\rfloor+1}^{m}
  \frac{n+kd-km}{m+n+kd-(k+1)i}{m+n+kd-(k+1)i\choose m-i} {ki-kd+i\choose i},\\
=&\binom{m+n}{m} \qquad \qquad \qquad \quad \hfill{(\mbox{by the Hagen-Rothe identity})}\\
&-\sum\limits_{i=\lfloor\frac{kd-1}{k+1}\rfloor+1}^{m}
  \frac{n+kd-km}{m+n+kd-(k+1)i}{m+n+kd-(k+1)i\choose m-i} {ki-kd+i\choose i},
\end{align*}
which is just \eqref{eq-NCp9} by collecting some terms.

\vspace{.2cm}
\noindent{\bf Acknowledgments.}
I would like to thank Arthur L.B. Yang and Bing-Qing Li for valuable comments, inspiring discussions and suggestions, and generous support on an earlier version. This work was supported by the National Natural Science Foundation of China Grant No. 71171018 and the Fundamental Research Funds for the Central Universities Grant FRF-BR-12-008.

\end{document}